%% file: euler-bernoulli.tex
\documentclass{article}

\usepackage{a4wide}

\usepackage{graphics,graphicx}
\usepackage{amsmath,amsthm}
\usepackage{amsfonts,mathrsfs,amssymb}

\newcommand{\Loneloc}{\ensuremath{L_{\text{loc}}^{1}}}
\newcommand{\wt}{\widetilde}
\def\pg{\mathhexbox278}
 
\input mymakros.tex

\title{Distributional solution concepts for the Euler-Bernoulli beam equation
with discontinuous coefficients\thanks{Supported by Ministry of Science of
Serbia, project 144016, and the Austrian Science Fund (FWF) START program Y237 on
'Nonlinear distributional geometry'}}

\author{\emph{G\"unther H\"ormann}\\
Fakult\"{a}t f\"{u}r Mathematik der Universit\"{a}t Wien\\
Nordbergstra{\ss}e 15, A-1090 Wien, Austria\\
\ \\
and\\
\ \\
\emph{Ljubica Oparnica}\\
Institute of Mathematics, Serbian Academy of Science\\
Kneza Mihajla 35, 11000 Belgrade, Serbia}

\date{April 2007}

\begin{document}

\maketitle

\begin{abstract}
We study existence and uniqueness of distributional solutions $w$ to the
ordinary
differential equation $\frac{d^{2}}{dx^{2}} \Big(  a(x)  \cdot \frac
{d^{2}w(x)}{dx^{2}}\Big) + P(x) \frac{d^{2}w(x) }{dx^{2}} = g(x)$ with
discontinuous
coefficients and right-hand side. For example, if $a$ and $w$ are non-smooth
the product
$a\cdot w''$ has no obvious
meaning. When interpreted on the most general level of the hierarchy of
distributional
products discussed in \cite[Chapter II]{O:92}, it turns out that existence
of a solution
$w$ forces it to be at least continuously differentiable. Curiously, the
choice of the
distributional product concept is thus incompatible with the possibility of
having a
discontinuous displacement function as a solution. We also give conditions
for unique
solvability.\\[2mm]
Key words: ordinary differential equations with discontinuous coefficients,
distributional solutions, multiplication of distributions.\\
AMS 2000 Subject Classification: 46F10;34A36.
\end{abstract}

\section{Equation of the Euler-Bernoulli beam}

We consider an Euler Bernoulli rod
under a distributed transversal force $g$ and axial force $P$.
The differential equation of equilibrium for the displacement
$w$ is given in \cite{Atanackovic:97} in the form
\begin{equation}
  \frac{d^{2}}{dx^{2}}\Big( EI \frac{d^{2}w(x)}{dx^{2}} \Big)
       +P \frac{d^{2} w(x)}{dx^{2}}=g(x)
         \qquad x\in\left[0,l\right]. \label{1.2}%
\end{equation}
Here, $E$ is the modulus of elasticity, $I$ is the moment of
inertia, and $l$ is the length of the rod. In our analysis we will
allow for nonconstant, $x$-depended, even discontinuous
coefficients $I$ and $P$. When there is a discontinuity in $I$ at
some point $x_{0}$ the rod can be considered to consist of two
different, but connected, parts, i.e., $EI(x) =
EI_{1}+H(x-x_{0})(EI_{2}-EI_{1})$ where $I_1 \not= I_2$ are the
corresponding moments of inertia respectively and $H$ denotes the
Heaviside function.

Equation (\ref{1.2}) has been studied in
\cite{YS:01}, where the authors discuss possible jump discontinuities at
$x_{0}$ in the
displacement
\begin{equation}
  \Delta := w(x_{0}+) - w(x_{0}-)   \label{1.6}%
\end{equation}%
(as well as in the rotation $\theta := w'(x_{0}+)-w'(x_{0}-)$), where a
suffix $+$ or
$-$ in the function argument denotes the limit from the right or left.
A solution Ansatz of the form $w=w_{1}+H\left(  x-x_{0}\right)  \left(
w_{2}-w_{1}\right)$ is then used, where $w_1$ and $w_2$ solve the equation
to the left
and to the right of $x_0$. In course of justifying this could be called a
solution $\Delta$ was being forced to vanish in order to avoid ill-defined
products involving a Dirac delta.

Here, we investigate the corresponding mathematical issues left
open: first, we analyze the possibility to give a meaning to the
notion of 'distributional solution' in the context of the
distributional product hierarchy described in \cite[Chapter
II]{O:92} (see also the Appendix for a brief review); second, we
show that indeed $\Delta$ necessarily has to vanish then, which is
consistent with the calculations in \cite{YS:01}; more precisely,
if $w$ were to have a jump discontinuity then the model product $[a \cdot w'']$,
which is the most general in the distributional product hierarchy, can not exist.
Thus, in order to allow for solutions with
jump discontinuities in the displacement one is forced to go beyond
intrinsic distributional products and use, e.g.\ algebras of
generalized functions (cf.\ \cite{GKOS:01,O:92}). For example, there has been
active research on such issues for hyperbolic partial differential equations with
discontinuous coefficients, where in certain cases non-existence of
distributional solutions
has been proved (cf.\ \cite{Hoermann:03,HdH:01,HdH:01c,LO:91,O:88,O:89}).

For notational simplification and structural clarity we put $A=EI_{1}$,
$B=EI_{2}$
(hence $A\neq B$) and
\begin{equation}
 a(x) = A + (B-A) H(x-x_{0}) =A H( x_{0}-x)
+B H(x-x_{0}). \label{1.12}
\end{equation}
Then the governing differential equation with boundary conditions  for the
Euler-Bernoulli rod with jump discontinuities in the bending read
\begin{equation}
  \frac{d^{2}}{dx^{2}}
        \Big( a(x) \cdot\frac{d^{2}w(x)}{dx^{2}} \Big) + P(x)
\frac{d^{2}w(x)}{dx^{2}} =
                g(x),    \qquad  x \in[0,1] \label{1.3}
\end{equation}
\begin{align}
w(0)  &  =0;\hspace{0.5cm}w(l)=0;\nonumber\\
\frac{d^{2}w}{dx^{2}}(0)  &  =0;\hspace{0.5cm}\frac{d^{2}w}{dx^{2}}(l)=0.
\label{1.4}
\end{align}
Mechanically, a global condition of equilibrium is expressed by equality of
the bending
moments
\begin{equation}
 E I_1 w''(x_0 -) = E I_2 w''(x_0 +). \label{1.5}
\end{equation}

We may use the substitution $u = w''$ to lower the order of equation
(\ref{1.3})\ and boundary conditions (\ref{1.4})
\begin{equation}
\frac{d^{2}}{dx^{2}}\Big( a(x) \cdot u(x) \Big)  + P(x)\, u(x)  =
  g(x)  \qquad x \in[0,1] \label{1.1}
\end{equation}
with
\begin{equation}
 u(0)=u(1)=0. \label{1.11}
\end{equation}
\begin{remark} Note that the above substitution is equivalent to imposing the
additional boundary problem $\frac{d^{2}}{dx^{2}}w(x)=u(x)$ with  $w(0)=w(1)=0$,
which is uniquely solvabile once $u$ is determined by (\ref{1.1}-\ref{1.11}). In
the sequel we will thus only consider $u$.
\end{remark}

  In equation (\ref{1.1}) the
product of the distributions $a$ and $u$ arises. There are several concepts
of partialy defined
products in the space of distributions. In the current paper we use the
 so-called model product (cf.\ \cite{O:92}) to give a meaning to the
differential
equation.

\begin{remark} {\bf (Comparison with $L^2$-operator theory.)} The above boundary
value problem (\ref{1.1}-\ref{1.11}) can as
well be investigated in the classical functional analytic context of unbounded
operators on $L^2([0,1])$. Singularities of the coefficient functions
then have a significant influence on choices for an appropriate
domain. In course of the current paper, we follow an intrinsic distribution
theoretic view, which allows for a wider class of solutions, right-hand sides in
the differential equation, as well as variations in the solution concept itself.

To illustrate the situation in an unbounded operator approach we briefly sketch
the constructions for the case where $a$ is given by (\ref{1.12}) and $P$ is a
real constant. It is natural to implement the boundary conditions (\ref{1.11})
into the domain of the operator. Furthermore, we have to specify the meaning of
the formal expression $(a u)''$. Note that requiring that $u$ belongs to the
Sobolev space $H^2(]0,1[)$ makes $u''$ well-defined in $L^2([0,1])$ and gives
sense to the boundary conditions $u(0) = u(1) = 0$. Observe that under these
hypotheses $a  u = A u_- + B u_+$, where $u_-$ (resp.\ $u_+$) vanishes to the
right (resp.\ left) of $x_0$ and is continuously differentiable on the left
(resp.\ right) up to $x_0$. Thus, by Schwartz' formula (\cite[Chapitre II,
\pg 2]{Schwartz:66}), we have $(a u)'' = a u'' + (B u_+'(x_0)-A u_-'(x_0))
\cdot \de_{x_0} + (B u_+(x_0) - A u_-(x_0)) \cdot \de_{x_0}'$, which is in $L^2$
only for $u$ such that the coefficients of $\de_{x_0}$ and $\de_{x_0}'$ vanish.

Therefore, we define the operator $T u := a \cdot u''$ with domain
$$
   D(T) := \{ u \in H^2(]0,1[) : u(0) = u(1) = 0, B u_-(x_0) - A u_+(x_0) = 0,
              B u_+'(x_0)-A u_-'(x_0) = 0 \}.
$$
It is straightforward to check that $T$ is
symmetric, i.e., $D(T) \subseteq D(T^*)$ and $T^* \mid_{D(T)} = T$, where $T^*$
denotes the adjoint of $T$. In fact, one
can prove that $T$ is self-adjoint along the following lines: Let $\vphi\in
\Cinfc(]0,1[) \cap D(T)$ and $v\in D(T^*)$; interpreting $L^2$-inner products
$\inp{.}{.}$ in terms of distributional actions $\dis{.\,}{.}$ and
vice versa we obtain
$\dis{\ovl{T^*v}}{\vphi} = \inp{\vphi}{T^* v}=\inp{T\vphi}{v} = \inp{a \vphi''}{v}
= \dis{\ovl{(a v)''}}{\vphi}$, which implies that $(a v)'' \in L^2$, forcing
that $v$ belongs to $H^2$ and satisfies the conditions appearing in $D(T)$ at
$x_0$. Furthermore, integration by parts is then applicable with $u\in D(T)$
yielding $\inp{T u}{v} = \inp{u}{a \cdot v''} + B u'(1) v(1) - A u'(0) v(0)$;
since $u \mapsto \inp{T u}{v}$ has to be a continuous linear
functional with respect to the $L^2$-norm $v(0)$ and $v(1)$ have to vanish. Hence
$v$ is in $D(T)$ and $T^* = T$.

We observe that the original differential operator in Equation (\ref{1.1}) is of
the form $T + P I$, where $I$ denotes the identity operator. Therefore, questions
concerning uniqueness and existence of solutions to (\ref{1.1}-\ref{1.11}) when
$g\in L^2$ directly relate to spectral properties of $T$. One can view
corresponding results obtained in Section 3 below in this context.
\end{remark}

\section{Solution concept based on the model product}

We analyze the properties of a distributional solution $u$ to
problem (\ref{1.1}-\ref{1.11}) in detail when the product $a \cdot
u$ is interpreted as a 'model product'. Throughout this and the
following two sections we focus on regularity issues stemming from
the highest order terms in the equation. Therefore we make the
assumption that

 \begin{center} $P$ is constant. \end{center}
We will remove this assumption and generalize our results in a final section
allowing
for jump discontinuities in $P$ as well.

\begin{definition}\label{d.1}
  Let $\dot{\D}'([0,1]) := \{ v\in\D'(\R); \supp v \subseteq [0,1]  \}$. We
call $u
  \in \dot{\D}'([0,1])$ a \emph{solution} to (\ref{1.1}-\ref{1.11}) if the
  following holds:

 \begin{description}
  \item{(A1)} The model product $[a \cdot u]$ of $u$ and $a$ (defined as in
  \cite{O:92}, see also the Appendix) exists in $\D'(\R)$

  \item{(A2)} The equation
   \begin{equation}
    (  [a \cdot u])'' + P u = g \label{2.1}
   \end{equation}
   holds in $\mathcal{D}^{\prime}( \mathbb{R})$.
 \end{description}
\end{definition}

\begin{remark}
 \begin{trivlist}
 \item{(i)} The boundary conditions (\ref{1.11}) are implemented into the
definition of
 the space of prospective solutions $\dot{\D}'([0,1])$ in the following
sense: if $u \in
 \dot{\D}'([0,1])$ happens to be a continuous function then $u(0) = u(1) =
0$.
 \item{(ii)} Note that (A1) is equivalent to the existence of the model
products $[
 H_{-}\cdot u] $ and $[ H_{+}\cdot u]$ where $H_{-}(x) = H(x_{0}-x)$ and
$H_{+} =
 H(x-x_{0})$.
 \end{trivlist}
\end{remark}

\begin{lemma}\label{T2}
\begin{trivlist}
\item{(i)} Let $u \in \dot{\D}'([0,1])$ satisfy (A1-2) then
$[H_{-}\cdot u]$ and $[ H_{+}\cdot u]$ belong to $\dot{\D}'([0,1])$.
\item{(ii)} $[ H_{-}\cdot\delta_{x_{0}}^{(k)}]$ and
$[H_{+}\cdot\delta_{x_{0}}^{(k)}]$
exist if and only if $k=0$, in which case we have
$[H_{-} \cdot\delta_{x_{0}}] = -\frac{\delta_{x_{0}}}{2}$,
$[ H_{+}\cdot\delta_{x_{0}} ] = \frac{\delta_{x_{0}}}{2}$. (Cf.\ similar
investigations in \cite[Lemma 4]{HdH:01})
\end{trivlist}
\end{lemma}

\begin{proof} Let $\vphi\in\D(\R)$ with $\int \vphi = 1$ and $\vphi_{\eps}(x) :=
\vphi(x/ \eps) / \eps$ be a model delta net (cf.\ \cite{O:92},(7.9)).

(i) By definition
$[H_{\pm}\cdot u] = \lim_{\eps\to 0}(H_{\pm} \ast \vphi_{\eps})
\cdot (u \ast \vphi_{\eps})$.
Let $w_{\eps} := (H_{\pm} \ast \vphi_{\eps}) \cdot (u \ast \vphi_{\eps})$ and
$\psi\in\D(\R)$ with $\supp \psi \cap [0,1] = \emptyset$ then
$\dis {w_{\eps}}{\psi} = \dis {H_{\pm} \ast \vphi_{\eps}}{(u \ast
\vphi_{\eps})\cdot\psi}$.
Since $\supp(u \ast \vphi_{\eps})\subseteq [0,1]+ \supp\vphi_{\eps} \subseteq
[-d_{\eps},1+d_{\eps}]$ for some  $d_{\eps} \to 0$ (as $\eps \to 0$) we have
$(u \ast \vphi_{\eps})\cdot\psi = 0$ and thus $\dis {w_{\eps}}{\varphi} =
0$.

(ii)  For any $\psi \in\D(\R)$
\begin{align*}
  \dis{ [ H_{-}\cdot\delta_{x_{0}}^{(k)} ] }{\psi} & =  \lim_{\eps\to 0}
       \dis{ ( H(x_{0} - x) \ast
         \vphi_{\eps})(\delta_{x_{0}}^{(k)} \ast\vphi_{\eps})}{\psi} \\
  & = \lim_{\eps\to 0} \frac{1}{\eps^{k+1}}
     \int\limits_{\R} \int\limits_{\frac{x-x_{0}}{\eps}}^{\infty}\vphi(t)
      \vphi^{(k)}(\frac{x-x_{0}} { \eps}) \psi(x)\, dt dx\\
  & = \lim_{\eps\to 0} \frac{1}{\eps^{k}}  \int\limits_{\R}
   \int\limits_{z}^{\infty}\vphi(t) \vphi^{(k)}(z) \psi(\eps z + x_{0}) \,
dt dz
\end{align*}
cannot be convergent for all $\psi$ as $\eps\to 0$ if $k \not= 0$. In case
$k=0$
we obtain the formula $[H_{-}\cdot\delta_{x_{0}}] = - \delta_{x_{0}}/2$ by
dominated convergence and the fact that $\int_{\R}\int_{z}^{\infty}\vphi(t)
\vphi(z)\, dz dt = - 1/2$.
The proof for $[ H_{+}\cdot\delta_{x_{0}} ]$ is similar.
\end{proof}

 \begin{theorem}\label{T1}
  Let $u \in \dot{\D}'([0,1])$ be a solution in the sense of Definition
\ref{d.1}.
  Then $u$ is a locally integrable function.
 \end{theorem}

 \begin{proof}
  \emph{Step 1:} Putting $\widetilde{u}_{-} = u \mid_{(0,x_{0})}$ and
   $\widetilde{u}_{+}=u \mid_{(x_{0},1)}$ yields
  \begin{align}
   A\widetilde{u}_{-}'' + P\widetilde{u}_{-} &= g \mid_{(0,x_{0})}
\label{2.2} \\
   B\widetilde{u}_{+}'' + P\widetilde{u}_{+} &= g \mid_{(x_{0},1)}.
\label{2.3}
  \end{align}
 Solving these two differential equations with constant coefficients we get
  \begin{equation}\label{2.4}
   \widetilde{u}_{-}=\widetilde{u}_{-h}+\widetilde{u}_{-p}\;
   \qquad\text{and}\qquad
   \widetilde{u}_{+} =\widetilde{u}_{+h}+\widetilde{u}_{+p},
  \end{equation}
 where
  \[
   \widetilde{u}_{-h}(x) = C_{1} e^{\sqrt{-P/A}x} +
       C_{2}e^{-\sqrt{-P/A}x}, \quad
   \widetilde{u}_{+h}(x) = D_{1}e^{\sqrt{-P/B}x} +
       D_{2}e^{-\sqrt{-P/B}x}
  \]
 and
 \[
 \widetilde{u}_{-p}(x) = \frac{1}{2\sqrt{-P/A}}(\int_{0}^{x}g(\tau)
     e^{\sqrt{-P/A}(x-\tau)}d\tau-
     \int_{0}^{x}g(\tau)e^{-\sqrt{-P/A}(x-\tau)}d\tau)
,%\label{p.r.}
\]
with a similar formula for $\widetilde{u}_{+p}(x)$ replacing $P/A$ by $P/B$
and
integration limits from $1$ to $x$.
Here, $\widetilde{u}_{-h}$, $\widetilde{u}_{+h}$ are smooth and
$\widetilde{u}_{-p}$, $\widetilde{u}_{+p}$ are absolutely continuous.
Therefore
$\widetilde{u}_{-}$ and $\widetilde{u}_{+}$ are absolutely continuous
functions on open subintervals $(0,x_{0})$ and $(x_{0},1)$ respectively.
Also, by
explicit formula, we see that $\widetilde{u}_{-}(x_{0}-) :=
\lim_{x\to x_{0}-}\widetilde{u}_{-}(x)$ and
$\widetilde{u}_{+}(x_{0}+)  := \lim_{x\to x_{0}+}\widetilde{u}_{+}(x)$
exist.

 \emph{Step 2:}  Define $\widetilde{u}\in \Loneloc(\R) $ by
  \begin{equation}\label{2.5}
   \widetilde{u}(x)  =
   \begin{cases}  0 & -\infty<x\leq0 \\
   \widetilde{u}_{-}(x) & 0<x<x_{0} \\
   \widetilde{u}_{+} (x) & x_{0}<x<1 \\
   0 & 1\leq x<\infty
   \end{cases}.
  \end{equation}
 We have that $\widetilde{u}\in\dot{\D}'([0,1])$ and
 $(u-\widetilde{u})\mid_{\R\backslash\{x_{0}\}} = 0$.
 Therefore $\supp(u-\widetilde{u}) =\{x_{0}\}$, which implies that
 \begin{equation}
 u=\widetilde{u}+\sum_{k=0}^{N}c_{k}\delta_{x_{0}}^{(k)},\qquad
c_{k}\in\C,N\in\N_{0}. \label{2.6}
 \end{equation}
 By Lemma \ref{T2} and Assumption (A1) $N=0$ in (\ref{2.6}). Hence
 \begin{equation}
  u=\widetilde{u}+c_{0}\delta_{x_{0}}. \label{2.8}
 \end{equation}

\emph{Step 3:} By Assumption (A2) we now obtain
 \begin{equation}
  ( u\cdot a)''=g-Pu=g-P\widetilde{u}-cP\delta_{x_{0}}, \label{2.7}
 \end{equation}
 where $g-P\widetilde{u}\in \Loneloc(\mathbb{R}) $.
 Let $w$ be a primitive function for $g-P\widetilde{u}$.
 Then $w-c_{0}PH(x-x_{0})$ is one for $(u\cdot a)'$.
 Therefore
  \[
   u\cdot a=W-c_{0}P(x-x_{0})_{+},
  \]
 where $x_{+}$ denote kink function, i.e.
   $x_{+}=\left\{\begin{array}
              [c]{c}
               x\\
               0
              \end{array}
              \begin{array}
               [c]{c}
               x>0\\
               x\leq0
              \end{array}
           \right.$
 and $W$ is primitive function for $w$.
 Since $ W\in C^{1} $ and the kink function is absolutely continuous we have
that $u\cdot a$ is
 absolutely continuous. But then (\ref{2.8}) and (\ref{1.12}) imply that
   \[
    u\cdot a=\widetilde{u}\cdot a+\frac{c_{0}}{2}\delta_{x_{0}}(B-A)
   \]
 which is absolutely continuous if and only if $c_{0}=0.$
 This in turn yields  $u=\widetilde{u}$ and therefore $u$ is locally
integrable
 as $\widetilde{u}$ is.
\end{proof}

\section{Existence and uniqueness of an $L^1([0,1])$-solution}

As we have seen in the previous section, a distributional solution
in the sense of Definition \ref{d.1} necessarily is a locally
integrable function. In this case, we can interpret the product
$u\cdot a$ as a duality product (cf.\ \cite{O:92} or the Appendix).
We analyze this situation more closely.

\begin{proposition}\label{M1}
 If $u\in L^1([0,1])$ is a solution to (\ref{1.1}-\ref{1.11}) then $u\in
C^{1}([0,1]\backslash\{x_{0}\})$
 and $u$ has a jump at $x=x_{0}$.
\end{proposition}

\begin{proof}
If $u\in L^1([0,1])$ then the differential equation (\ref{1.1}) yields
 $(u\cdot a)''\in L^1([0,1])$, hence $(u\cdot a)'\in
C_{\text{abs}}([0,1])$
 and thus $u\cdot a\in C^{1}([0,1])$.
 Therefore we also have that
 $u\cdot a\mid_{[0,x_{0})}=Au\mid_{[0,x_{0})} \in C^{1}([0,x_{0}))$ in turn
$u\in C^{1}([0,x_{0}))$.
 Similarly, $u\in C^{1} ((x_{0},1])$.
 Furthermore $\lim_{x\to x_{0}-}u\cdot a(x)=A\cdot u(x_{0})$ and therefore
 $u(x_{0}-)=\lim_{x\to x_{0}-}u(x)$ exists.
 Similarly for $u(x_{0}+) =\lim_{x\to x_{0}+}u(x)$. But $a\cdot u$ is
continuous, so that
 $\lim_{x\to x_{0}-}u\cdot a(x) =\lim_{x\to x_{0}+}u\cdot a(x)$ and thus
   \begin{equation}\label{3.1}
     Au(x_{0}-) = Bu(x_{0}+),
   \end{equation}
 which implies the global equilibrium condition (\ref{1.5}).
 If $u\in C([0,1])$ then (\ref{3.1}) implies $A=B$,
 which contradicts the assumption $I_{1}\neq I_{2}$.
 This means that $u$ has to be discontinuous at $x_{0}$.
\end{proof}

\begin{remark}\label{AC2}
 As a matter of fact we have $u'\in C_{\text{abs}}([0,1]\backslash
\{x_{0}\})$.
 Indeed, since $(u\cdot a)'\in C_{\text{abs}}([0,1])$ reasoning as above
 we obtain that $u'$ is absolutely continuous off $x_{0}$, so that
 $u'(x_{0}-)$, $u'(x_{0}+)$ exist and obtain
  \begin{equation}\label{3.11}
   Au'(x_{0}-) = Bu'(x_{0}+).
  \end{equation}
\end{remark}

 Now we are in a position to construct a solution to (\ref{1.1}-\ref{1.11}).

 \begin{lemma}\label{P_lemma} For any choice of $A > 0$, $B > 0$, and $0< x_0
<1$ there exists a strictly increasing sequence $(P_l)_{l\in\N}$
  of positive real numbers $P_l$ such that the following holds:
\begin{enumerate}
\item[(i)] If $P < 0$
\item[or]
\item[(ii)] if $P >0$ and $P \not= P_l$ for all $l\in\N$
\end{enumerate}
then there is a unique solution to (\ref{2.2}) and (\ref{2.3})
with $\widetilde{u}_{-}(0)=0 $ and $\widetilde{u}_{+}(1)=0 $, which satisfies the
stability  conditions
    \begin{align}
      A\wt{u}_{-}(x_0)=B\wt{u}_{+}(x_0)  \label{sk} \\
      A\wt{u}'_{-}(x_0)=B\wt{u}'_{+}(x_0). \label{sk1}
    \end{align}
 \end{lemma}

 \begin{remark} \label{excases}
  In case $P = P_{l}$ for some $l\in\N$ the solution is not unique or even may
fail to exist. Investigation of these cases seems possible in a direct way
without requiring further analytical tools.
 \end{remark}

 \begin{proof}  Any solutions to (\ref{2.2}) and (\ref{2.3})  are  given by
(\ref{2.4}).

 \emph{Case $P < 0$:} The solution formulae (\ref{2.4}), adapted to the
boundary
  conditions at $0$ and $1$, give
    $$\wt{u}_{-}(x) = 2C_{1}\sinh \sqrt{-P/A}\,x + \wt{u}_{-p}(x) $$
  and
    $$\wt{u}_{+}(x) = 2D_{1}e^{\sqrt{-P/B}x}\sinh \sqrt{-P/B}\,(x-1) +
\wt{u}_{+p}(x), $$
  where
   \[
     \wt{u}_{-p}(x) = \frac{1}{\sqrt{-P/A}}\int_{0}^{x}g(\tau)
       \sinh \sqrt{-P/A}(x-\tau) \, d\tau
   \]
 and similarly to $\wt{u}_{+p}(x)$  (replacing $P/A$ with $P/B$ and
integration limits
 from $1$ to $x$). The stability conditions (\ref{sk}-\ref{sk1}) are
equivalent to the
 linear system $Hy=z$ with
    \begin{equation}
      H := \left[\begin{array}
               [c]{cc}
            2A\sinh\sqrt{-P/A}x_{0}
            & -2Be^{\sqrt{-P/B}}\sinh\sqrt{-P/B}(x_{0}-1) \\
            2A\sqrt{-P/A}\cosh\sqrt{-P/A}x_{0}
 &-2B\sqrt{-P/B}e^{\sqrt{-P/B}}\cosh\sqrt{-P/B}(x_{0}-1)
        \end{array}\right]
     \end{equation}
  and
   \[
    y := \left[\begin{array}[c]{c}
        C_{1}\\
        D_{1}
      \end{array} \right],\qquad
    z := \left[\begin{array}[c]{c}
        Bu_{+p}(x_{0}) - Au_{-p}(x_{0})\\
        Bu_{+p}'(x_{0}) - Au_{-p}'(x_{0})
      \end{array}\right].
   \]
  Further we have
   \begin{align*}
      \det H  & =4AB e^{\sqrt{-P/B}}
                 \left(-\sqrt{-P/B}\sinh
                  \sqrt{-P/A}x_{0}\cosh\sqrt{-P/B}(x_{0}-1)
                 \right. \nonumber\\
               & \left. +\sqrt{-P/A}\sinh\sqrt{-P/B}
                  (x_{0}-1)\cosh\sqrt{-P/A}x_{0}
                  \right).  \label{5.1}
     \end{align*}
  Since   $$\sinh\sqrt{-A/P}x_{0}\cosh\sqrt{-B/P}(x_{0}-1) > 0 $$
  and     $$\sinh\sqrt{-B/P}(x_{0}-1)\cosh\sqrt{-A/P}x_{0} < 0 $$
  we have $ \det H<0 $
  hence unique solvability of the above linear system.

\emph{Case $P >0$:}
   If $P>0$ then the solutions $\wt{u}_{-}$, $\wt{u}_{+}$ involve $\sin$ and
$\cos$
   (instead of $\sinh$ and $\cosh$) and the determinant of the corresponding
linear  system $H y = z$ reads
    \begin{align*}
      \det H  & =4AB e^{\sqrt{-P/B}}
                 \left(-\sqrt{P/B}\sin
                  \sqrt{P/A}x_{0}\cos\sqrt{P/B}(x_{0}-1)
                 \right. \nonumber\\
               & \left. +\sqrt{P/A}\sin\sqrt{P/B}
                  (x_{0}-1)\cos\sqrt{P/A}x_{0}
                  \right).
     \end{align*}
  When $\det H \not= 0$ we have the same situation as in the case $P<0$.
Observe that the
  set $Z_0$ of values for $P$ such that any cosine factor occurring in the
above  determinant  vanishes is at most countable. Apart from these values, to
find
   $P>0$ for which $\det H = 0 $ is equivalent to solving
   $$ h(s) := \tan(s) + \nu \mu \tan(\mu s) = 0, $$
  where $s=\sqrt{P/A} x_0 > 0$, $ \mu = \sqrt{A/B} (1/x_0 - 1) > 0$,
  and $\nu = x_0 / (1 - x_0) > 0$. One observes that there is a countable
discrete set
  of singularities of $h$, at which the limits from the left and right are
$+\infty$
  and $-\infty$ respectively. Since $h$ is continuous otherwise, there is a
countable
  set $Z_1$ of (positive) zeroes. To summarize, the union $Z_0 \cup Z_1$
makes up a
  sequence $(P_l)_{l\in\N}$ with the required property.
 \end{proof}

 \begin{remark} We point out that the above proof of Lemma \ref{P_lemma}
does not give the
  minimum set of values $P_l$ to be removed. In fact, only those elements in
$Z_0$ have to
  occur in $(P_l)$ which make both cosine factors vanish. Note that the
latter can only
  happen, when $\sqrt{B/A} x_0 / (1 - x_0)$ is a rational number of the form
  $(2l+1)/(2k+1)$ with integers $k$, $l$.
 \end{remark}

 \begin{theorem}\label{T3}
Let  $P$ satisfy one of the conditions (i)-(ii) in Lemma \ref{P_lemma}.   Let
$\wt{u}_{-}$ and $\wt{u}_{+}$ be the solutions to (\ref{2.2})  and (\ref{2.3})
obtained in Lemma \ref{P_lemma} and  define
    \begin{equation}\label{3.2}
    u_{-}(x)=\left\{\begin{array}
                 [c]{c}
                 \widetilde{u}_{-}(x)  \\
                 0
               \end{array}
               \begin{array}
                 [c]{c}
                 x\in\left[  0,x_{0}\right] \\
                 x\in\left[  x_{0},1\right]
               \end{array},\right. \qquad
   u_{+}(x)=\left\{\begin{array}
                 [c]{c}
                 0\\
                 \widetilde{u}_{+}(x)
              \end{array}
              \begin{array}
                 [c]{c}
                 x\in [0,x_{0}] \\
                 x\in [x_{0},1]
              \end{array} \right. .
   \end{equation}
  Then
     \begin{equation}\label{3.3}
      u(x) = u_{-}(x) + u_{+}(x)
     \end{equation}
   is the unique solution to (\ref{1.1}-\ref{1.11}) in the sense of
Definition \ref{d.1},
  belongs to $L^1([0,1])$ and satisfies the boundary conditions in the
classical
  sense.
 \end{theorem}

  \begin{proof}
   Since $a(x) = AH(x_{0}-x) + BH(x-x_{0})$ we have
     $$
       (u\cdot a)(x) = Au_{-}(x) + Bu_{+}(x)
     $$
  which we will differentiate twice.
  Recall (\cite[Chapitre II, \pg 2]{Schwartz:66}) that if a function $f$ is in
$C_{\text{abs}}([0,1]\backslash\{x_{0}\})$, such
  that $\lim_{x\to x_{0}-} f(x) = f(x_{0}-)$ and $\lim_{x\to x_{0}+}f(x) =
f(x_{0}+)$
  exist, then the distributional derivative $\frac{d}{dx}f$ satisfies
   \begin{equation}\label{3.5}
     \frac{d}{dx}f(x) = f'(x) + \big(f(x_{0}+)-
f(x_{0}-)\big)\cdot\delta_{x_{0}},
   \end{equation}
  where $f'$ denotes the (class of) function(s) in $L^1([0,1])$ equal to
the pointwise
  derivative of $f$ almost  everywhere in  $[0,1]\setminus\{ x_{0} \}$.

  Therefore
   \[
    \frac{d}{dx}(u\cdot a) =
    A u_{-}'+Bu_{+}'+\big(B u(x_{0}+)- Au(x_{0}-)\big)\cdot\delta_{x_{0}}
   \]
  and
   \begin{multline}\label{3.6}
     \frac{d^{2}}{dx^{2}}(u\cdot a) = Au_{-}'' + Bu_{+}''+
                \big(B u'(x_{0}+) - A u'(x_{0}-)\big)\cdot\delta_{x_{0}}\\
                 + \big(B u(x_{0}+) - Au(x_{0}-)\big)\cdot\delta_{x_{0}}'.
   \end{multline}
 By construction we have that
  $$
  Au_{-}''(x)=\left\{\begin{array}
                   [c]{c}
                   (-Pu+g)(x)\\
                   0
                  \end{array}
                  \begin{array}
                    [c]{c}
                    x\in[0,x_{0}] \\
                    x\in[x_{0},1]
                  \end{array}\right., \;
  B u_{+}''(x) =\left\{\begin{array}
                   [c]{c}
                   0\\
                   (-Pu+g)(x)
                  \end{array}
                  \begin{array}
                   [c]{c}
                   x\in[0,x_{0}] \\
                   x\in[x_{0},1]
                  \end{array}\right..$$
   Thus (\ref{3.1}) and (\ref{3.11}) imply that
    \[
     \frac{d^{2}}{dx^{2}}(u\cdot a) = -Pu + g.
    \]
   Note that $u$ is continuous near the boundaries $x=0$ and $x=1$, thus
   the conditions (\ref{1.11}) follow by construction.

For uniqueness, we first observe that any solution $u$ has to be in
$\Loneloc$ by
  Theorem \ref{T1}. Furthermore, due to Proposition \ref{M1} it also has to
satisfy the
  stability conditions (\ref{3.1}-\ref{3.11}). Hence Lemma \ref{P_lemma}
implies
  uniqueness.
  \end{proof}

 \section{Generalization to discontinuous axial force}

 We extended the analysis of the previous section to investigate solvability
of the same type of differential equation
  \[
   \frac{d^{2}}{dx^{2}}[a(x)\cdot u(x)]+ Pu(x) = g(x)
  \]
 with boundary condition $u(0)  =  u(1)  = 0$,  where the force $P$ now is a jump
function of the form
  \[
   P = P_{1}+ H(x-x_{0})(P_{1}-P_{2})
  \]
  with real numbers $P_1$, $P_2$.

 As with constant $P$ we obtain that any solution to the
 differential equation (in a sense similar to Definition \ref{d.1})
necessarily is a locally integrable function and
 continuously differentiable off $x_{0}$ with a jump at $x=x_0$.

 \begin{remark}\label{r2}
   Note that condition (A1) in Definition \ref{d.1} implies that the
model product $[P \cdot u]$ exists. Therefore we will now require $u$ to be a
solution to the differential equation in the sense of Definition \ref{d.1} with
(\ref{2.1}) replaced by
   \[
   ([a \cdot u])'' + [P \cdot u] = g.
   \]
  \end{remark}

 \begin{theorem}
  (i) Let $u\in  \dot{\D}'([0,1])$ be a solution in the sense of Remark
\ref{r2}.
  Then $u$ is a locally integrable function.

  (ii) If $u\in \Loneloc$ is a solution then $u\in
C^{1}([0,1]\backslash\{x_{0}\})$
  and $u$ has a jump discontinuity at $x=x_{0}$. Furthermore, equations
(\ref{3.1}) and
(\ref{3.11}) hold.
 \end{theorem}

 \begin{proof}
  \emph{Step 1:} As in the proof of Theorem \ref{T1} we can set
  $\widetilde{u}_{-}=u\left|_{( 0,x_{0})  }\right. $
  and $\widetilde{u}_{+}=u\left|_{( x_{0},1)}\right.$.
  Then we have
   \begin{align}
    A\widetilde{u}_{-}'' + P_{1}\widetilde{u}_{-} &
=g\left|_{(0,x_{0})}\right.\label{6.6}\\
    B\widetilde{u}_{+}'' + P_{2}\widetilde{u}_{+} &
=g\left|_{(x_{0},1)}\right.\label{6.7}
   \end{align}
  Solving these two differential equations with constant coefficients we get
   \begin{equation}\label{6.3}
    \wt{u}_{-}=\wt{u}_{-h} + \wt{u}_{-p}
    \hspace{0.5cm}\text{and}\hspace{0.5cm}
    \widetilde{u}_{+}=\widetilde{u}_{+h}+\widetilde{u}_{+p}\
   \end{equation}
  where
   \begin{equation*}
    \wt{u}_{-h}(x) =C_{1}e^{\sqrt{-P_{1}/A}x} +
    C_{2}e^{-\sqrt{-P_{1}/A}x},\;
    \wt{u}_{+h}(x) =D_{1}e^{\sqrt{-P_{2}/B}x} +
D_{2}e^{-\sqrt{-P_{2}/B}x}
   \end{equation*}
  and
   \[
   \wt{u}_{-p}(x) = \frac{1}{2\sqrt{-P_{1}/A}}
    \big( \int_{0}^{x}g(\tau)e^{\sqrt{-P_{1}/A}(x-\tau)}d\tau
    -\int_{0}^{x}g(\tau)e^{-\sqrt {-P_{1}/A}(x-\tau)}d\tau\big),
   \]
  with a similar formula for $\widetilde{u}_{+p}(x)$ replacing $P_{1}/A$ by
$P_{2}/B$ and
  integration limits from $1$ to $x$. Again,
$\widetilde{u}_{-h},\widetilde{u}_{+h}$ are smooth and
  $ \wt{u}_{-p}$,$\wt{u}_{+p}$ are absolutely continuous.
  Therefore $\wt{u}_{-}$ and $\wt{u}_{+}$
  are absolutely continuous on open subintervals $(0,x_0)$ and $(x_0,1)$.
  Also, there exist $\wt{u}_{-}(x_{0-}) = \lim_{x\rightarrow
x_{0-}}\wt{u}_{-}(x)$
  and $\wt{u}_{+}(x_{0+}) = \lim_{x\rightarrow x_{0+}}\wt{u}_{+}(x)$.

  \emph{Step 2:} Precisely as in Step 2 of the proof of Theorem \ref{T1} we
obtain
   \begin{equation}\label{6.2}
    u=\widetilde{u}+c\delta_{x_{0}}.
   \end{equation}

  \emph{Step 3:} Equation (\ref{6.2}) and $P=P_{1}+ H(x-x_{0})(
P_{1}-P_{2})$
  leads to $ Pu=P_{1}\wt{u}_{-}+ P_{2}\wt{u}_{+} +
\frac{c}{2}\delta_{x_{0}}(P_{2}-P_{1})$.
  Since $P_{1}\wt{u}_{-}+P_{2}\wt{u}_{+}\in \Loneloc(\R)$
  we have that $w=\int(g-P_{1}\wt{u}_{-}-P_{2}\wt{u}_{+})$ is
  absolutely continuous and its primitive function $W$ is $C^{1}$.
  The differential equation $( u\cdot a)'' = g - Pu $
  implies
   \[
    u\cdot a=W-\frac{c}{2}\left(  P_{2}-P_{1}\right)  \left(  x-x_{0}\right)
 _{+},
   \]
  which is absolutely continuous. The same arguments as in Theorem \ref{T1}
yield that $c=0$ and hence $u=\widetilde{u}$ is locally integrable. This proves
(i).

  For part (ii) we may reason as in the proof of Proposition \ref{M1} we get
that
  second part of theorem is valid.
 \end{proof}

 The construction of a solution rests on the following lemma which corresponds to
Lemma \ref{P_lemma}.

 \begin{lemma}\label{P1_lemma} For any choice of $A > 0$, $B > 0$, and $0< x_0
<1$ there exist one-dimensional submanifolds $\cal{M}$ and $\cal{N}$ of $\R^2$
such that the following holds:
\begin{enumerate}
 \item[(i)] If $P_1 < 0$ and $P_2 < 0$
 \item[or]
 \item[(ii)] if $P_1 > 0$, $P_2 > 0$ and $(P_1,P_2) \not\in \cal{M}$
 \item[or]
 \item[(iii)] if $P_1 > 0$, $P_2 < 0$, and $(P_1,P_2) \not\in \cal{N}$
 \end{enumerate}
 then there is a unique solution to (\ref{6.6}) and (\ref{6.7}) with
  $\widetilde{u}_{-}(0)=0 $ and $\widetilde{u}_{+}(1)=0 $, which satisfies the
stability conditions
    \begin{align}
      A\wt{u}_{-}(x_0)=B\wt{u}_{+}(x_0)  \label{ask} \\
      A\wt{u}'_{-}(x_0)=B\wt{u}'_{+}(x_0). \label{ask1}
    \end{align}
 \end{lemma}

 \begin{remark}
Similarly as in Remark \ref{excases} for the cases where $(P_1,P_2)$ belongs to
$\cal{M}$ or $\cal{N}$ the solution is not unique or may
fail to exist, explicit investigation of which could be carried out
along the lines of the following proof.
 \end{remark}

 \begin{proof}
 As in the proof of Lemma \ref{P_lemma} any solution to (\ref{6.6}) and (\ref{6.7}) is given
by (\ref{6.3}).

  \emph{In case $P_1 < 0$ and $P_2 < 0$} the solution
formulae (\ref{6.3}) with boundary conditions at $0$ and $1$, and stability conditions (\ref{ask})
and (\ref{ask1}) lead to a linear ($2$x$2$) system $ Hy=z $ with $y$ and $z$ as in Lemma
\ref{P_lemma} and
   \begin{align*}
      \det H  & =4AB e^{\sqrt{-P_2/B}}
                 \left(-\sqrt{-P_2/B}\sinh
                  \sqrt{-P_1/A}x_{0}\cosh\sqrt{-P_2/B}(x_{0}-1)
                 \right. \nonumber\\
               & \left. +\sqrt{-P_1/A}\sinh\sqrt{-P_2/B}
                  (x_{0}-1)\cosh\sqrt{-P_1/A}x_{0}
                  \right) < 0.
     \end{align*}
   Therefore we have a unique solution.
\vspace{0.2cm}

\noindent
\parbox[c]{.6\linewidth}{\ \hphantom{bla}\emph{Case $P_1 > 0$, $P_2 > 0$}:
$\det H$ now reads
   \begin{align*}
      \det H  & =4AB e^{\sqrt{P_2/B}}
                 \left(-\sqrt{P_2/B}\sin
                  \sqrt{P_1/A}x_{0}\cos\sqrt{P_2/B}(x_{0}-1)
                 \right. \nonumber\\
               & \left. + \sqrt{P_1/A}\sin\sqrt{P_2/B}
                  (x_{0}-1)\cos\sqrt{P_1/A}x_{0}
                  \right).
   \end{align*}
Whenever this is nonzero we have a unique solution.
To see where it vanishes let $ s =
   \sqrt{P_1/A}x_{0} $, $ t = \sqrt{P_2/B}(x_{0}-1)$, $\nu = x_0/(1-x_0)$ and
analyze the function
  \begin{equation*}
    f(s,t)=\nu t \sin s\cos t + s\sin t\cos s.
  \end{equation*}
By direct inspection one deduces that $\grad f$ is
nonzero when $f=0$
 which yields that the zero set ${\cal M} '=\{ (s,t)\in{\R^2}; f(s,t)=0 \}$ is a
one-dimensional
 submanifold  of $\R^2$.
 We set ${\cal M} = \{ (P_1,P_2)\in{\R^2}; (s,t)\in\cal M '  \} $.}
\parbox{.35\linewidth}{\includegraphics[width=5cm,height=5cm]{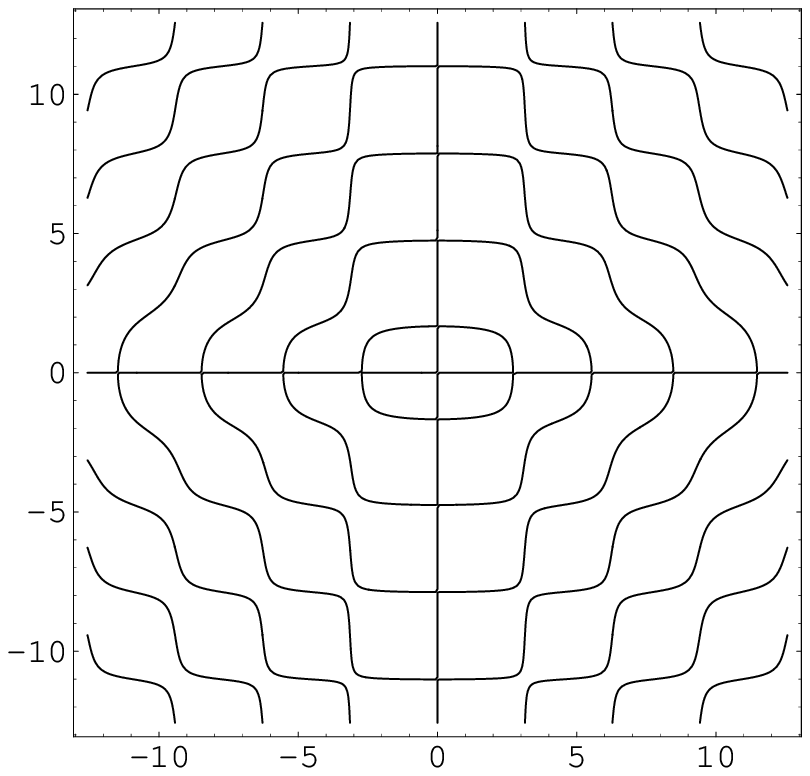}\\
{\small Figure: $\cal{M}'$ is a union of infinitely
many closed concentric curves (plot for the case $\nu =
6$).}}
\vspace{0.1cm}

   \emph{Case  $P_1 > 0$, $P_2 < 0$}: we have
   \begin{align*}
      \det H  & =4AB e^{\sqrt{P_2/B}}
                 \left(\sqrt{-P_2/B}\sinh
                  \sqrt{-P_1/A}x_{0}\cos\sqrt{P_2/B}(x_{0}-1)
                 \right. \nonumber\\
               & \left. + \sqrt{P_1/A}\sin\sqrt{P_2/B}
                  (x_{0}-1)\cosh\sqrt{P_1/A}x_{0}
                  \right).
   \end{align*}
 and as above one can show that the zero set is a one-dimensional
 submanifold $\cal N$ of $\R^2$. In the complement $\det H \not= 0$, solution exists and is unique.
 \end{proof}

 \begin{theorem} Let $P_1$ and $P_2$ satisfy one of the conditions (i)-(iii) in
Lemma \ref{P1_lemma}.
  Let $\wt{u}_{-}(x)$ and $\wt{u}_{+}(x)$ be the solutions to (\ref{6.6}) and
(\ref{6.7}) obtained in Lemma \ref{P1_lemma} and define
   \begin{equation} \label{6.5}
    u_{-}(x)=\left\{\begin{array}
              [c]{c}
              \wt{u}_{-}(x)\\
               0
              \end{array}
              \begin{array}
              [c]{c}
              x\in [0,x_{0}] \\
              x\in [x_{0},1]
              \end{array}\right.
    \hspace{0.5cm}
    u_{+}(x)=\left\{\begin{array}
              [c]{c}
              0\\
              \widetilde{u}_{+}(x)
              \end{array}
              \begin{array}
              [c]{c}
              x\in [0,x_{0}] \\
              x\in [x_{0},1]
             \end{array}\right. .
   \end{equation}
Then
   \begin{equation}\label{6.4}
    u(x) = u_{-}(x) + u_{+}(x)
   \end{equation}
 is a solution to (\ref{1.1}-\ref{1.11}) with
 $P=P_{1}+H(x-x_{0})(P_{1}-P_{2})$.
 \end{theorem}

 \begin{proof}
  Proceeding as in the proof of Theorem \ref{T3} we arrive at (\ref{3.6}). By
construction we have
   $$Au_{-}'' = \left\{\begin{array}
    [c]{c}
    -P_{1}u+g\\
     0
    \end{array}
     \begin{array}
      [c]{c}
      x\in\left[  0,x_{0}\right]  \\
      x\in\left[  x_{0},1\right]
     \end{array}
    \right.
    \qquad \text{and}\qquad
    Bu_{+}''=\left\{\begin{array}
     [c]{c}
     0\\
     -P_{2}u+g
     \end{array}
     \begin{array}
      [c]{c}
      x\in[0,x_{0}]\\
      x\in[x_{0},1]
      \end{array}
      \right.   .
$$
Employing (\ref{ask}) and (\ref{ask1}) we get
$\frac{d^{2}}{dx^{2}}(u\cdot a) = -(P_{1}+(P_{2}-P_{1})H(x-x_{0}))u + g$.
\end{proof}

\section{Approximation by regularization}

In this section we investigate the possibility to approximate the solution
to (\ref{1.1}), (\ref{1.11}) using some regularization of the coefficient
$a(x)$. Throughout this section we will assume that $P$ is constant and such
that the solution to (\ref{1.1}), (\ref{1.11}) is unique.

 Let $Qu:=(au)''+ Pu$ and consider the equation $Qu = g\in L^1([0,1])$.
 Suppose that $a_{\eps}$ is a smooth regularization of the jump disconuity in $a$ such that
 $a_\eps \in C^2([0,1])$ with $\lim_{\eps\to 0} a_{\eps}^{(j)} = a^{(j)}$ uniformly on compact subsets of $[0,1]\setminus\{x_0\}$ for derivative orders $j = 0, 1,2$. Let $Q_{\eps}u := (au_{\eps})''+ Pu$.

 \begin{proposition} Let $g\in L^1([0,1])$.\\[1mm] 
 (i) If $u_{\eps} \in \dot{\D}'([0,1])$ denotes the solution to 
$$ 
   Q_{\eps}u = g, \qquad u(0)=0, u(1)=0,
$$
then $u_{\eps}$ belongs to the space $AC^2([0,1])$ of continuously diifferentiable functions whose derivates are loccaly integrable.\\[1mm]
(ii) Let $u$ be the solution to $Qu = g$. Then $u_{\eps}\to u$ uniformly on compact subsets of $[0,1]\setminus\{x_0\}$.
\end{proposition}
 \begin{proof}
 (i): 
  Note that $a_{\eps} Q_{\eps} u_{\eps} = a_{\eps}g$ is 
  equivalent to the Sturm
  Liouville problem $L u_{\eps} := (pu_{\eps}')' + q u_{\eps} = a_{\eps} g$
   with $p := a_{\eps}^2$ and
  $ q := a_{\eps}(a_{\eps}'' + P)$. If $Q_{\eps} u_{\eps} = a_{\eps} g$ then
  $(a_{\eps}^2 u_{\eps}')' = a_{\eps} g - a_{\eps}(a_{\eps}'' + P)\in L^1([0,1]) $ 
  since $a_{\eps} g\in L^1([0,1])$ and $a_{\eps}$ is $C^2$. Therefore
  $a_{\eps}^2 u_{\eps}'$ is absolutely continuous. Since $a_{\eps}$ is bounded from below away from zero this
  implies that $u_{\eps}'$ is absolutely continuous as well. Thus 
  $u_{\eps}\in AC^2([0,1])$.\\[1mm]
(ii): Let $v_{\eps}:= u_{\eps}-u$. Then we obtain
 \begin{equation} \label{Qeps_eq_feps}
  Q_{\eps}v_{\eps} = ([(a_{\eps}-a)\cdot u])'' : = f_{\eps}, 
  \qquad v_{\eps}(0)=v_{\eps}(1)=0.
\end{equation}
By assumption and using that $a_{\eps} - a$ is $C^2$ off $x_0$  we have that $f_{\eps}\to 0$ in $\Loneloc([0,1]\setminus\{x_0\})$.
Integrating twice in (\ref{Qeps_eq_feps}) gives
 \begin{equation}
  (a_{\eps}v_{\eps})(x) = \int_0^x (x-y)f_{\eps}(y)dy + P \int_0^x
  (x-y)v_{\eps} dy + a_{\eps}(0)v_{\eps}(0) +
  (a_{\eps}v_{\eps})'(0)\, x.
 \end{equation}
 
 Let $K$ be a compact subset of $[0,x_0)$. 
 Since $u$ and $u_{\eps}$ are $AC^2$ on $K$ (as noted in Remark \ref{AC2} and the prooof above) and both belong to $\dot{\D}'([0,1])$ we have that $v_{\eps}(0) = v'_{\eps}(0)=0$. Furthermore, by  
 $|(a_{\eps}v_{\eps})(x)|\geq A |v_{\eps}| > 0$ we obtain
 $$
  |v_{\eps}(x)|\leq \frac{1+P}{A}\int_0^x |x-y||f_{\eps}(y)|dy + \frac{1+P}{A}\int_0^x
  |x-y||v_{\eps}(y)|dy.
 $$
Applying Gronwall's inequality we get
 $$
  |v_{\eps}(x)| \leq \frac{1+P}{A}\, e^{\frac{(1+P) x}{A}}\int_0^x |x-y||f_{\eps}(y)|dy
  \leq\frac{1+P}{A}\, e^{\frac{1+P}{A}}\, \|f_{\eps}\|_{L^1(K)}  .
 $$
As noted above, we have $f_{\eps} \to 0$ in $\Loneloc([0,1]\setminus\{x_0\})$, thus $v_{\eps} \to 0$ uniformly on $K$. 
  
The reasoning in case of a compact subset contained in $(x_0,1]$ is similar. Since an arbitrary compact subset of $[0,1]\setminus\{x_0\}$ is the disjoint union of two compact subsets in either part of $[0,1]\setminus\{x_0\}$ the assertion is proved.
\end{proof}

We illustrate the convergence in an example: we put $A=1$, $B =2$, $P=1$, $x_0 = 1/2$,  
and design $a_{\eps}$ as the $C^2$ function, which is defined by a fifth order odd polynomial in $[-\eps,\eps]$ and equal to $a$ otherwise. As right-hand side we choose $g(x) = - \cos(11 x) / \sqrt{|x - 2/3|}$. Note that $g \in L^1([0,1]) \setminus L^2([0,1])$. The following plots show the regularized coefficients and corresponding solutions for parameter values $\eps = 1/10, 1/30, 1/100$. 

\vspace{0.5cm}

\parbox{.4\linewidth}{\includegraphics[width=6cm,height=4cm]{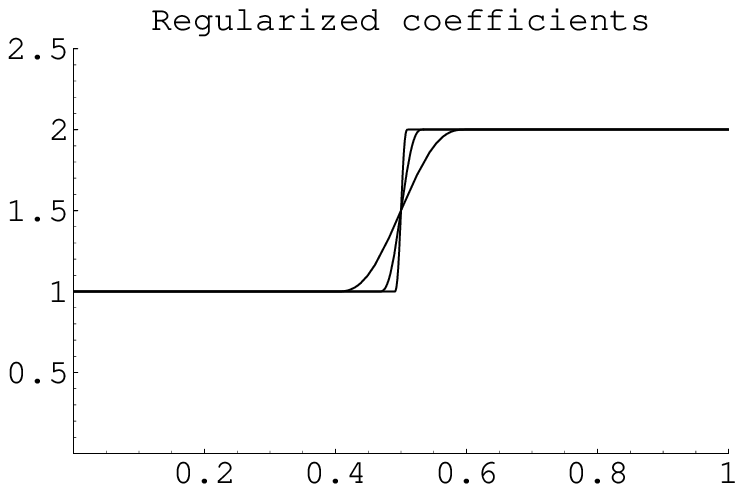}}
\hfill
\parbox{.5\linewidth}{\includegraphics[width=6cm,height=4cm]{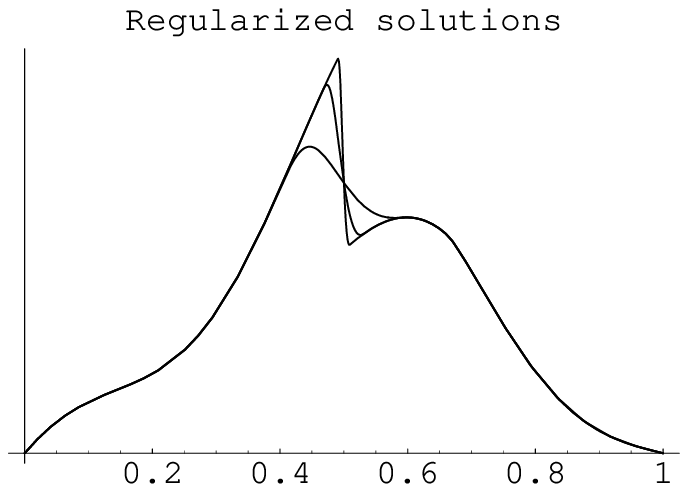}}

\vspace{0.1cm}

\appendix

\section{Appendix: Hierarchy of distributional products}

For convenience of the reader we briefly review the basic definitions of the
coherent distributional products described in
\cite[Chapter II]{O:92} in terms of a hierarchy.
All these products yield the classical multiplication when
restricted to smooth functions.

We use $\Omega$ to denote an open subset of $\R^n$ and $\FT{u}$ for the Fourier
transform of $u$.

The most elementary product in this context is
$\mathbf{\Cinf\cdot\D'}$, the product of a smooth function and a
distribution, defined as the adjoint of multiplication by a smooth
function in the test function space.

\paragraph{Disjoint singular support:}
Assume that $u$, $v$ are in $\D'(\Omega)$ with disjoint singular supports.
Then for any $x\in\Omega$ there is a neighborhood $\Omega_x$ and a function
$f_x\in\D(\Omega_x)$ such that either $f_x u$ or $f_x v$ is smooth. Then
in $\Omega_x$ the product of $u$ and $v$ can be defined in the sense of
$\Cinf\cdot\D'$ and by the localization properties of $\D'$
(cf.~\cite{Hoermander:V1},
subsect.~2.2) this consistently defines a distribution in $\Omega$.

\paragraph{Wave front set condition:}
Let $u\in\D'(\Omega)$
and $(x_0,\xi_0)\in T^*\Omega\setminus 0 := \{ (x,\xi) \mid x\in\Omega, \xi
\not= 0 \}$ (the cotangent bundle over $\Omega$ with the zero section
removed). $u$ is said to be \emph{microlocally regular} at $(x_0,\xi_0)$ if
there is $\vphi\in\D(\Omega)$, $\vphi(x_0) \not= 0$, and an open cone
$\Gamma$
with axial vector $\xi_0$ such that $\FT{\vphi u}$ is rapidly decreasing in
$\Gamma$. $\WF(u)$ is the closed subset of $T^*\Omega\setminus 0$ where $u$
is
not microlocally regular.

If $u$, $v$ $\in\D'(\Omega)$ their wave
front sets are said to be in \emph{favorable position} if $(x,\xi)\in\WF(u)$
implies
that $(x,-\xi)\not\in\WF(v)$. In this case the product of $u$ and $v$ can be
defined as the pullback of the tensor product $u\otimes v
\in\D'(\Omega\times\Omega)$ by the diagonal map
$\Omega\to\Omega\times\Omega$,
$x\mapsto(x,x)$ (cf.~\cite{Hoermander:V1}, Thm.~8.2.10).

\paragraph{Fourier product:}
Given two distributions $u$, $v$ $\in\D'(\Omega)$ we say that their Fourier
product exists if for every $x\in\Omega$ there is an open neighborhood
$\Omega_x$ and $f_x\in\D(\Omega)$, $f_x = 1$ on $\Omega_x$, such that the
${\S}'$-convolution of $\FT{f_x u}$ and $\FT{f_x v}$ exists. Locally near
$x$,
the product of $u$ and $v$ is then defined to be the inverse Fourier
transform
of $\FT{f_x u} * \FT{f_x v}$ (for a definition of $\S'$-convolvability see
\cite{O:92}, sect.~6).

\paragraph{Duality products:}
Let $X$ be a \emph{normal space of distributions}, that is $\D \subseteq X
\subseteq \D'$ and $\D$ is dense in $X$. Assume that the dual space $X'$ is
(equipped with a locally convex topology so that it becomes) normal as well
and that multiplication with a fixed element in $\D$ induces a continuous
linear map both from $X$ into $X$ and from $X'$ into $X'$.

For any normal space of distributions $Y$ denote by $Y_{\text{loc}}$ the set
of distributions $v\in\D'$ such that $\psi v \in Y$ for all $\psi\in\D$. If
$u\in(X')_{\text{loc}}$ and $v\in X_{\text{loc}}$ then the product of $u$
and
$v$ can be defined by
\[
  \dis{u\cdot v}{\psi} := \dis{\chi u}{\psi v}
\]
for $\psi\in\D$ and $\chi\in\D$ chosen arbitrarily with
$\chi = 1$ on $\supp(\psi)$. Note that in the above definition the left hand
side denotes a $(\D',\D)$ pairing while the right hand side uses
the pairing $(X',X)$.

\paragraph{Strict and model products:}
The basic idea is to regularize one or both factors by convolution, perform
the multiplication in the sense $\Cinf\cdot \D'$ or $\Cinf\cdot\Cinf$, and
try
to take the limit. The regularizing convolutions are carried out with two
principal types of mollifiers.

A net $(\rho^\eps)_{\eps > 0}$ in
$\D(\R^n)$ is called \emph{strict delta net} if
\begin{align}
\supp(\rho^\eps) \to \{0\} \quad \text{ as } \eps \to 0 \\
\int \rho^\eps(x) \, dx = 1 \quad \text{ for all } \eps > 0 \\
\int |\rho^\eps(x)|\, dx \quad \text{is bounded independently of $\eps$}.
\end{align}

A \emph{model delta net} is given by specifying $\vphi\in\D(\R^n)$ with
$\int\vphi(x)\, dx = 1$ and defining $(\vphi_\eps)_{\eps > 0}$ by $\vphi_\eps(x)
= \vphi(x/\eps)/\eps^n$.

Consider the following four possibilities to define a product of $u$ and
$v$:
\begin{eqnarray} \setcounter{equation}{1}
 u \cdot [v] &=& \lim\limits_{\eps\to 0} u (v*\rho^\eps) \\
 {[u]} \cdot v &=& \lim\limits_{\eps\to 0} (u*\rho^\eps) v \\
 {[u]} \cdot [v] &=& \lim\limits_{\eps\to 0} (u*\rho^\eps) (v*\sig^\eps) \\
 {[u \cdot v]} &=& \lim\limits_{\eps\to 0} (u*\rho^\eps) (v*\rho^\eps)
\end{eqnarray}
where the limit is required to exist in $\D'(R^n)$ and independent of
the choice of $(\rho^\eps)_{\eps > 0}$ and $(\sig^\eps)_{\eps > 0}$ in the
class of strict, resp.\ model, delta nets. This defines four types of so-called
\emph{strict}, resp.~\emph{model}, products. Since the definitions (1)-(3)
turn
out to be equivalent when using strict, resp.~model, delta nets
(cf.~\cite{O:92}, Thms.~7.2 and 7.11)
we distinguish only the following four products: strict product (1)-(3),
strict product (4), model product (1)-(3), and model product (4).

%\vspace{-0.4cm}

\noindent\parbox{0.35\textwidth}{\paragraph{Coherence properties:}
The various products satisfy coherence properties and can be brought into
the following hierarchy table. Here, an arrow indicates that a product
definition is contained and consistent with its successor in the graph. All
products shown generalize the multiplication $\Cinf \cdot \D'$.
}
\hphantom{Bla}
\parbox{0.6\textwidth}{
\input hier.tex
}
%\vspace{1cm}

\paragraph{Acknowledgement:} The authors are very grateful to Prof.\ Teodor
Atanackovic
for providing the mechanical model and valuable suggestions for mathematical
investigations. We also thank Simon Haller for critical discussions on
several details in the construction and regularity of the solution.

\bibliographystyle{abbrv}
\bibliography{buba}

\end{document}

%% file: mymakros.tex
\newcommand{\Cinf}{\ensuremath{\mathcal{C}^\infty}}
\newcommand{\Cinfc}{\ensuremath{\mathcal{C}^\infty_{\text{c}}}}
\newcommand{\D}{\ensuremath{{\cal D}}}
\renewcommand{\S}{\mathscr{S}}

\newcommand{\mb}[1]{\ensuremath{\mathbb{#1}}}
\newcommand{\N}{\mb{N}}

\newcommand{\R}{\mb{R}}
\newcommand{\C}{\mb{C}}

\newcommand{\WF}{\mathrm{WF}}

\newcommand{\grad}{\ensuremath{\mbox{\rm grad}\,}}

\newfont{\bl}{msbm10 scaled \magstep2}

\newtheorem{theorem}{Theorem}[section]
\newtheorem{lemma}[theorem]{Lemma}
\newtheorem{proposition}[theorem]{Proposition}
\newtheorem{definition}[theorem]{Definition}

\theoremstyle{definition}
\newtheorem{remark}[theorem]{Remark}

\newcommand{\beq}{\begin{equation}}
\newcommand{\eeq}{\end{equation}}
\renewenvironment{proof}[1][Proof]{\textbf{#1.} }{\ \rule{0.5em}{0.5em}}

\newcommand{\FT}[1]{\widehat{#1}}

\newcommand{\dis}[2]{\langle #1 , #2 \rangle}
\newcommand{\inp}[2]{\langle #1 | #2 \rangle}  %math mode
\newcommand{\notmid}{\mid\kern-0.5em\not\kern0.5em}

\newcommand{\de}{\delta}
\newcommand{\eps}{\varepsilon}

\newcommand{\vphi}{\varphi}

\newcommand{\sig}{\sigma}

\newcommand{\supp}{\mathop{\mathrm{supp}}}

\newcommand{\ovl}[1]{\overline{#1}}

%% file: hier.tex
\setlength{\unitlength}{1500sp}
\begingroup\makeatletter\ifx\SetFigFont\undefined%
\gdef\SetFigFont#1#2#3#4#5{%
  \reset@font\fontsize{6}{#2pt}%
  \fontfamily{#3}\fontseries{#4}\fontshape{#5}%
  \selectfont}%
\fi\endgroup%
\begin{picture}(9179,5761)(316,-5132)
\thicklines
\put(2971,-3526){\vector( 1, 0){765}}
\put(2971,-2176){\vector( 1, 0){765}}
\put(2971,-871){\vector( 1, 0){765}}
\put(5086,-1141){\vector( 0,-1){765}}
\put(5086,-2536){\vector( 0,-1){765}}
\put(5086,-3886){\vector( 0,-1){765}}
\put(4276,-1006){\makebox(0,0)[lb]{\smash{\SetFigFont{10}{12.0}{\rmdefault}{
\mddefault}{\updefault}WF favorable}}}
\put(4276,-2356){\makebox(0,0)[lb]{\smash{\SetFigFont{10}{12.0}{\rmdefault}{\mddefault}{\updefault}Fourier product}}}
\put(4276,-5056){\makebox(0,0)[lb]{\smash{\SetFigFont{10}{12.0}{\rmdefault}{\mddefault}{\updefault}model product (1)-(3)}}}
\put(4276,-3661){\makebox(0,0)[lb]{\smash{\SetFigFont{10}{12.0}{\rmdefault}{\mddefault}{\updefault}strict product (1)-(3)}}}
\put(9451,-3841){\vector( 0,-1){765}}
\put(8596,-3661){\makebox(0,0)[lb]{\smash{\SetFigFont{10}{12.0}{\rmdefault}{\mddefault}{\updefault}strict product (4)}}}
\put(8596,-5011){\makebox(0,0)[lb]{\smash{\SetFigFont{10}{12.0}{\rmdefault}{\mddefault}{\updefault}model product (4)}}}
\put(7471,-3571){\vector( 1, 0){765}}
\put(7471,-4966){\vector( 1, 0){765}}
%\put(1351,209){\vector( 0,-1){765}}
\put(1351,-1096){\vector( 0,-1){765}}
\put(1351,-2446){\vector( 0,-1){765}}
\put(451,-3616){\makebox(0,0)[lb]{\smash{\SetFigFont{10}{12.0}{\rmdefault}{\mddefault}{\updefault}$W^{m,p}_{\text{loc}}$-duality}}}
%\put(811,389){\makebox(0,0)[lb]{\smash{\SetFigFont{10}{12.0}{\rmdefault}{
%\mddefault}{\updefault}$\Cinf\cdot\D'$}}}
\put(541,-2221){\makebox(0,0)[lb]{\smash{\SetFigFont{10}{12.0}{\rmdefault}{\mddefault}{\updefault}$H^s_{\text{loc}}$-duality}}}
\put(316,-961){\makebox(0,0)[lb]{\smash{\SetFigFont{10}{12.0}{\rmdefault}{\mddefault}{\updefault}disjoint sing supp}}}
\end{picture}